\documentclass[12pt,draft]{amsart}

\makeatletter

\theoremstyle{plain}
\newtheorem{thm}{Theorem}[section]
\numberwithin{equation}{section} 
\numberwithin{figure}{section} 
\theoremstyle{plain}
\newtheorem*{thm*}{Theorem}
\theoremstyle{plain}
\theoremstyle{plain}
\newtheorem*{cor*}{Corollary}

\theoremstyle{plain}
\newtheorem{lem}[thm]{Lemma} 
\theoremstyle{plain}
\newtheorem{prop}[thm]{Proposition} 
\theoremstyle{definition}

\theoremstyle{remark}

\theoremstyle{remark}

\theoremstyle{remark}

\theoremstyle{remark}

\theoremstyle{definition}

\theoremstyle{remark}
  \newtheorem*{acknowledgement*}{Acknowledgement}

\theoremstyle{plain}

\theoremstyle{plain}
\theoremstyle{plain}
\theoremstyle{plain}
\theoremstyle{definition}

\theoremstyle{remark}

\theoremstyle{remark}

\theoremstyle{remark}

\theoremstyle{plain}

\newcommand{\p}{\varphi}

\usepackage{amssymb} \usepackage{amscd}

\makeatother

\begin{document}

\title{Inductive limits, unique traces and tracial rank zero}

\author{Nathanial P. Brown}

\address{Department of Mathematics, Penn State University, State
College, PA 16802}

\email{nbrown@math.psu.edu}

\thanks{Partially supported by DMS-0244807.}

\begin{abstract}
In the program to classify C$^*$-algebras, it is very important to find abstract conditions which are sufficient to imply that a given algebra has tracial rank zero, in the sense of Huaxin Lin.  Even in the presence of a unique trace, we show that the union of the known necessary conditions is not enough. 
\end{abstract}

\maketitle

\section{Introduction}

In \cite{lin:classify} Huaxin Lin made a breakthrough in Elliott's classification program: C$^*$-algebras of tracial rank zero are amenable to classification.  This remarkable theorem has since been applied in a variety of contexts, illustrating the (pleasantly!) surprising fact that checking the tracial-rank-zero axioms is possible in many \emph{concrete} examples.  

At the other end of the spectrum, it is natural and important to search for \emph{abstract} hypotheses which would imply that a particular class of algebras has tracial rank zero.  Evidently the scope of Huaxin's classification theorem would be substantially broadened by such a result.  

The obvious place to start, when looking for the `right' abstract hypotheses, would be the necessary ones: Every simple, unital, separable C$^*$-algebra of tracial rank zero enjoys the following properties: 
\begin{enumerate} 
\item[$\bullet$]  Real rank zero,  stable rank one and quasidiagonality \cite[Theorem 3.4]{lin:TAF};

\item[$\bullet$]  The Riesz interpolation property \cite[Theorem 6.11]{lin:TR}; 

\item[$\bullet$]  The fundamental (tracial) comparison property\footnote{This is, at least formally, stronger than Blackadar's original formulation which used quasitraces \cite{blackadar}.} -- i.e.\ if $p,q \in A$ are projections and $\tau(p) < \tau(q)$ for every tracial state on $A$ then $p$ is (Murray-von Neumann) equivalent to a subprojection of $q$ \cite[Theorem 6.8]{lin:TR}; 

\item[$\bullet$]  There exists an increasing sequence of residually finite dimensional subalgebras with dense union \cite[Theorem 3.8]{lin:AC}. 
\end{enumerate}

This is the essential list of known necessary conditions. (Other important properties follow from these, like weak unperforation and cancellation of projections, of course.)

It is known that even if $A$ is an exact C$^*$-algebra with all the properties above, it need not have tracial rank zero (cf.\ \cite[Theorem 6.2.4]{brown:invariantmeans}). The reason is von Neumann algebraic: every \emph{tracial} GNS representation of an algebra with tracial rank zero must be hyperfinite. (This is immediate from the definition of tracial rank zero, together with the old fact, due to Murray and von Neumann, that `locally' finite dimensional implies AFD \cite[Chapter 4]{MvN}.)  Indeed, the example constructed in \cite[Theorem 6.2.4]{brown:invariantmeans} has a non-hyperfinite II$_1$-factor representation, hence can't have tracial rank zero.  

On the other hand, it follows from \cite[Theorems 3.2.2 and 4.3.3]{brown:invariantmeans} that an exact quasidiagonal C$^*$-algebra \emph{with unique trace} must produce the hyperfinite II$_1$-factor in its GNS representation. Hence the obstruction vanishes in the unique trace case -- so long as $A$ is exact.  Thus it is still open, and exceedingly important to decide whether or not every exact algebra with the properties above, and possessing a unique tracial state, has tracial rank zero.  An affirmative answer would be a major breakthrough in the classification program; a counterexample would be devastating. 

Inspired by a question of Lin, we construct in this paper  the first example with all the properties above -- plus unique trace --  which doesn't have tracial rank zero.   The obstruction is again von Neumann algebraic, hence our example is not exact.  More precisely, the main result of this note is: 

\begin{thm} 
\label{thm:main} There exists a unital, separable, simple C$^*$-algebra $A$, containing a dense nest of RFD subalgebras (hence, is quasidiagonal), and with real rank zero, stable rank one, the fundamental (tracial) comparison property, Riesz interpolation, and a unique trace whose GNS representation yields a non-hyperfinite II$_1$-factor.  Thus $A$ does not have tracial rank zero.\footnote{Actually, its tracial rank is infinity -- i.e.\ it has no nice approximations at all.  Indeed, any `tracially nuclear' C$^*$-algebra will always produce hyperfinite tracial GNS representations.} 
\end{thm} 

The construction is very similar to that in \cite[Section 6.2]{brown:invariantmeans} and hence is heavily influenced by Dadarlat's seminal work on nonnuclear tracially AF algebras  \cite{dadarlat}. (We also reuse the main idea from \cite[Proposition 9.3]{brown:qdsurvey}.)

Being rather technical, we don't  feel that traditional exposition is the best way to convey the proof.  The next section outlines the main ingredients, highlighting crucial points without worrying about truth: i.e.\ we describe what we would \emph{like} to do, but don't explain why it's possible to do it.  Even in Section \ref{sec:works} we don't prove it's possible, we prove it works. (That is, \emph{if} one could carry out the procedure in Section \ref{sec:properties} then a C$^*$-algebra satisfying all the required hypotheses  exists.) In Section \ref{sec:details} we tidy up, explaining why Section \ref{sec:properties} is not a big hypothetical heap of rubbish.

\section{The Construction: Abstract Properties}
\label{sec:properties}

\subsection*{Data and Notation}
We begin with a description of the initial data. One needs a subalgebra $$E \subset \prod_{n \in \mathbb{N}} M_{k_n}\footnote{$M_{k_n}$ will denote the $k_n \times k_n$ complex matrices and $\prod_{n \in \mathbb{N}} M_{k_n}$ is the von Neuman algebra of bounded sequences.}$$ such that:
\begin{enumerate}
\item $\bigoplus M_{k_n} \triangleleft E$ and $E$ is separable, unital and has real rank zero (cf.\ \cite{brown-pedersen});\footnote{$\bigoplus M_{k_n}$ is the ideal of sequences tending to zero in norm.}

\item $E/(\bigoplus M_{k_n})$ has a unique tracial state $\tau_{\infty}$, and all matrix algebras over $E/(\bigoplus M_{k_n})$ have comparison with respect to their unique traces;

\item $\pi_{\tau_{\infty}}(E)^{\prime\prime}$ is \emph{not} hyperfinite, where $\pi_{\tau_{\infty}}$
denotes the GNS representation.
\end{enumerate}

The central projections in $E$ (coming from $\bigoplus M_{k_n} $) will play an important role, so let's give them a name: $1_{k_n}$ will denote the unit of $0\oplus \cdots \oplus 0 \oplus M_{k_n} \oplus 0 \cdots$ (but keep in mind that this is just a central projection in $E$, not a unit). We also need the infinite rank complements.  That is, let $$P_s = 1_E - \big(1_{k_1} \oplus 1_{k_2} \oplus \cdots \oplus 1_{k_s}\big).$$ We then define $E_s = P_sE$ and our picture becomes $$E = M_{k_1} \oplus E_1 = M_{k_1}\oplus M_{k_2} \oplus E_2 = M_{k_1}\oplus M_{k_2} \oplus M_{k_3} \oplus  E_3 = \cdots.$$

\subsection*{General structure and properties}
With the data in hand, we will then define natural numbers $l(s)$ and projections $r_s \in M_{l(s)}$ such that
\begin{enumerate}
\item[(4)]  $\lim_{s\to \infty} \mathrm{tr}(r_s \otimes r_{s-1} \otimes \cdots \otimes r_1) > 0$, \end{enumerate}
where $\mathrm{tr}$ is the unique tracial state on a matrix algebra.

Finally, we will construct an inductive system $$E = E_0 \to M_{l(1)}\otimes E_1 \to M_{l(2)}\otimes M_{l(1)} \otimes E_2 \to M_{l(3)} \otimes M_{l(2)}\otimes M_{l(1)} \otimes E_3 \to \cdots $$ by defining $*$-homomorphisms $$\p_{i}\colon M_{l(i)}\otimes \cdots \otimes M_{l(1)} \otimes E_i \to M_{l(i+1)}\otimes \cdots \otimes M_{l(1)} \otimes E_{i+1}$$ with all of the following properties:
\begin{enumerate}
\item[(5)] Each $\p_i$ is unital and injective;

\item[(6)] For every $0\neq x \in M_{l(i)}\otimes \cdots \otimes M_{l(1)} \otimes E_i$ there exists $j > i$ and a projection $w \in M_{l(j)}\otimes \cdots \otimes M_{l(1)} \otimes E_j$ such that $$0\neq w\p_{j,i}(x) \in M_{l(j)}\otimes \cdots \otimes M_{l(1)} \otimes 1_{E_j},$$ where $\p_{j,i} = \p_{j-1}\circ \p_{j-2} \circ \cdots \circ \p_{i}$;

\item[(7)] The projections $r_j\otimes \cdots \otimes r_{i+1}\otimes 1_{l(i)} \otimes \cdots \otimes 1_{l(1)}\otimes 1_{E_j}$ commute with $\p_{j,i}(M_{l(i)}\otimes \cdots \otimes M_{l(1)} \otimes E_i)$;

\item[(8)] And finally, for all $x \in M_{l(i)}\otimes \cdots \otimes M_{l(1)} \otimes E_i$ we have  $$(r_j\otimes \cdots \otimes r_{i+1}\otimes 1_{l(i)} \otimes \cdots  \otimes 1_{l(1)}\otimes 1_{E_j})\p_{j,i}(x) = r_j\otimes \cdots \otimes r_{i+1}\otimes \pi_{j,i}(x),$$ where $\pi_{j,i} \colon M_{l(i)}\otimes \cdots \otimes M_{l(1)} \otimes E_i \to M_{l(i)}\otimes \cdots \otimes M_{l(1)} \otimes E_j$ is the identity map on the matrices tensored with the projection map $E_i \to E_j$, $x \mapsto P_j x$.\end{enumerate}

\section{Why it works}
\label{sec:works}

Let $A$ denote the inductive limit of our hypothetical sequence $$E_0 \to M_{l(1)}\otimes E_1 \to M_{l(2)}\otimes M_{l(1)} \otimes E_2 \to M_{l(3)} \otimes M_{l(2)}\otimes M_{l(1)} \otimes E_3 \to \cdots.$$  Let $\Phi_i\colon M_{l(i)}\otimes \cdots \otimes M_{l(1)} \otimes E_i \to A$ denote the canonical $*$-homomorphisms.

Since each $E_i$ is residually finite dimensional and has real rank zero, it follows immediately from properties (1) and (5) that $A$ is a unital, separable C$^*$-algebra with real rank zero and containing a nested sequence of residually finite dimensional subalgebras with dense union.  That leaves 6 things to check: simplicity, stable rank one, Riesz interpolation, unique trace, not tracial rank zero and the comparison property. 

\subsection*{Simplicity}

This argument is well known and follows from (6): any ideal in $A$ would intersect some $\Phi_i(M_{l(i)}\otimes \cdots \otimes M_{l(1)} \otimes E_i)$ which, after pushing out to $M_{l(j)}\otimes \cdots \otimes M_{l(1)} \otimes E_j$ and multiplying by $w$, implies the ideal intersects a  \emph{unital} matrix subalgebra of $A$ -- thus contains the unit of $A$ (since matrix algebras are simple).

\subsection*{Stable rank one and the Riesz property}

These follow from the fact that our construction yields an approximately divisible C$^*$-algebra (cf.\ \cite[Theorem 1.4 and Corollary 3.15]{bkr}). If this isn't obvious -- actually, it isn't obvious until we describe the connecting maps explicitly.  

No problem.  If you can't wait, just replace $A$ with $A\otimes \mathcal{U}$, where $\mathcal{U}$ is the CAR algebra, say, and note that $A\otimes \mathcal{U}$ satisfies all the desired properties and obviously is approximately divisible.

\subsection*{Unique trace}

This is the meat.  It boils down to (2), which is the key to the following lemma.

\begin{lem}
\label{lem:cuttrace}
Let $\tau$ be any tracial state on $A$ and $\pi_{\tau}\colon A\to B(H)$ the corresponding GNS representation. For each $i \in \mathbb{N}$ there exists a projection $R_i \in \pi_{\tau}(A)^{\prime\prime}$ such that
$$\tau(R_i \pi_{\tau}(\Phi_i(x))) = \tau(R_i) (\mathrm{tr}\otimes \tau_{\infty})(x),$$ for all $x
\in M_{l(i)}\otimes \cdots \otimes M_{l(1)} \otimes E_i$, and $\tau(R_i) \to 1$ as $i \to \infty$.
\end{lem}

\begin{proof}  For each $i$ we define $R_i$ to be the weak limit of the decreasing sequence of projections $$\pi_{\tau}(\Phi_j(r_j\otimes \cdots \otimes r_{i+1} \otimes 1_{l(i)}\otimes \cdots \otimes 1_{l(1)}\otimes 1_{E_j})),$$ as $j \to \infty$.  (Decreasing is not automatic, it follows from condition (8).)

These projections tend to 1 in trace because condition (4) -- and uniqueness of traces on matrix algebras -- ensures that $$\lim_{i \to \infty} \bigg(\lim_{j \to \infty} \tau(\Phi_j (r_j\otimes \cdots \otimes r_{i+1} \otimes 1_{l(i)}\otimes \cdots \otimes 1_{l(1)}\otimes 1_{E_j}))\bigg) = 1.$$

By continuity and relation (8) we have
\begin{eqnarray*}
\tau(R_i \pi_{\tau}(\Phi_i(x))) &=&  \lim_{j\to \infty}  \tau\bigg(\Phi_j\big((r_j\otimes \cdots \otimes r_{i+1} \otimes 1_{l(i)}\otimes \cdots \otimes 1_{l(1)}\otimes 1_{E_j}) \p_{j,i}(x)\big)\bigg)\\
&=& \lim_{j\to \infty} \tau\bigg(\Phi_j( r_j\otimes \cdots \otimes r_{i+1}\otimes \pi_{j,i}(x))\bigg)\\
&=& \lim_{j\to \infty} \mathrm{tr}(r_j\otimes \cdots \otimes r_{i+1}) \tau\bigg(\Phi_j(1_{l(j)}\otimes \cdots \otimes 1_{l(i+1)} \otimes \pi_{j,i}(x))\bigg)\\
\end{eqnarray*}
for all $x \in M_{l(i)}\otimes \cdots \otimes M_{l(1)} \otimes E_i$. However, since $$M_{l(i)}\otimes \cdots M_{l(1)} \otimes  (E/(\bigoplus M_{k_n}))$$ has a unique trace -- namely $\mathrm{tr}\otimes \tau_{\infty}$, thanks to (2) -- it follows that $$\lim_{j\to \infty} \tau\bigg(\Phi_j(1_{l(j)}\otimes \cdots \otimes 1_{l(i+1)} \otimes \pi_{j,i}(x))\bigg) = (\mathrm{tr}\otimes \tau_{\infty})(x)$$ and thus the limit of the product is the product of the limits, as desired.
\end{proof}

\begin{prop} $A$ has a unique tracial state.
\end{prop}

\begin{proof} $A$ must have a tracial state because it is the norm closure of an increasing union of subalgebras which have traces. (Any weak-$*$ limit of traces must be a trace.)

Now suppose $A$ has two traces, $\tau_1$ and $\tau_2$.  Let $R^{(1)}_i$ and $R^{(2)}_i$ be the projections from the previous lemma.  Then for every $x \in M_{l(i)}\otimes \cdots \otimes M_{l(1)} \otimes E_i$ of norm one, we use the fact that  $\tau(a)= \tau(aP) + \tau(aP^{\perp})$ ($P$ a projection) to deduce
\begin{eqnarray*}
|\tau_1(\Phi_i(x)) - \tau_2(\Phi_i(x))| & \leq & |\tau_1(R_i^{(1)}\pi_{\tau_1}(\Phi_i(x))) - \tau_2(R_i^{(2)}\pi_{\tau_2}(\Phi_i(x)))| + \epsilon_i\\
&=& |(\mathrm{tr}\otimes \tau_{\infty})(x)| |\tau_1(R_i^{(1)}) - \tau_2(R_i^{(2)})| + \epsilon_i,
\end{eqnarray*}
where $\epsilon_i = 2 - (\tau_1(R_i^{(1)}) + \tau_2(R_i^{(2)})).$ Evidently this implies $\tau_1 = \tau_2$.
\end{proof}

\subsection*{Non-hyperfinite GNS representation}

Let $\tau$ denote the unique trace on $A$. Lemma \ref{lem:cuttrace} implies that the von Neumann algebra generated by
the subalgebra (cf.\ (7)) $$R_i\pi_{\tau}(\Phi_i(M_{l(i)}\otimes \cdots \otimes M_{l(1)} \otimes E_i))$$ is isomorphic
to $$M_{l(i)}\otimes \cdots \otimes M_{l(1)} \otimes \pi_{\tau_{\infty}}(E)^{\prime\prime}.$$  Since the latter is not
hyperfinite (condition (3)), we deduce that $A$ can't have tracial rank zero (subalgebras of finite, hyperfinite von Neumann algebras must also be hyperfinite,
thanks to Connes' remarkable theorem \cite{connes:classification}).

\subsection*{The fundamental (tracial) comparison property}

This is the potatoes: not particularly interesting, just a necessary, rather bland part of the meal.  

\begin{lem} 
The algebra $E$ has the fundamental (tracial) comparison property. 
\end{lem} 

\begin{proof} Let $\tau_n$ be the tracial state on $E$ gotten by composing the coordinate projection $x \mapsto 1_{k_n}x$ with the trace on $M_{k_n}$.  Slightly abusing notation, we let $\tau_{\infty}$ also denote the trace on $E$ coming from $E/(\bigoplus M_{k_j})$.  

Now assume $p,q \in E$ are projections such that $\tau(p) < \tau(q)$, for all tracial states $\tau$ on $E$.  In particular this holds for $\tau_{\infty}$ -- and $E/(\bigoplus M_{k_j})$ has comparison by assumption -- so we can find a large integer $N$ such that $pP_N$ is equivalent to a subprojection of $qP_N$. (The details here are standard and left to the reader.  The key point is that partial isometries in $E/(\bigoplus M_{k_j})$ can be lifted to partial isometries in $E$.)  To fix the first $N-1$ coordinates, we use the traces $\tau_1, \ldots, \tau_{N-1}$ and the fact that $E$ contains the ideal $\bigoplus M_{k_j}$. 
\end{proof} 

Evidently the lemma above can be generalized to matrices over $E$.  Hence, $A$ is an inductive limit of algebras which enjoy the fundamental comparison property with respect to traces.   

\begin{prop} $A$ has the fundamental (tracial) comparison property. 
\end{prop} 

\begin{proof}  This argument is well-known, so we only sketch the main ingredients. 

Let $p, q\in A$ be projections such that $\tau(p) < \tau(q)$ for all tracial states $\tau$ on $A$.  Assume $A_n \subset A_{n+1}$ are subalgebras with the (tracial) comparison property.  We may assume, after perturbing, that $p, q \in A_n$, for some large $n$.  

We claim that there exists $m > n$ such that $\gamma(p) < \gamma(q)$ for all tracial states $\gamma$ on $A_m$ (and this will evidently complete the proof).  Indeed, if not we can find traces $\gamma_m$ on $A_m$ such that $\gamma_m(p) \geq \gamma_m(q)$ for all $m$.  Passing to a subsequential limit, this implies the existence of a trace $\tau$ on $A$ such that $\tau(p) \geq \tau(q)$.  Contradiction. 
\end{proof}

\section{How to do it}
\label{sec:details}

Now comes the fun part.  Let's start with the data.

\subsection*{Existence of data}

Our requisite algebra exists because of the following theorem.

\begin{thm}  There exists a separable, unital MF algebra\footnote{By definition, this means a (separable, unital) subalgebra of the quotient $\frac{\prod M_{k_n}}{\bigoplus M_{k_n}}$ for some choice of natural numbers $k_n$.} $B$ with real rank zero, unique tracial state $\tau_{\infty}$, comparison with respect to $\tau_{\infty}$, and such that $\pi_{\tau_{\infty}}(B)^{\prime\prime}$ is not hyperfinite.
\end{thm}

Arranging all of these properties simultaneously is quite a deep fact.  Indeed, the reduced group C$^*$-algebra $C^*_r(\mathbb{F}_2)$ of a free group is MF by \cite{ht:ext}.  Hence, so is $C^*_r(\mathbb{F}_2)\otimes \mathcal{U}$, for any UHF algebra $\mathcal{U}$.  But this latter algebra has real rank zero and comparison with respect to its unique trace, by  work of R{\o}rdam (cf.\ \cite{rordam}, \cite{bkr}) and Haagerup \cite{haagerup}.  Evidently this does the trick, since free group factors are not hyperfinite. 

Since there exist $k_n$ such that $$B \subset \frac{\prod M_{k_n}}{\bigoplus M_{k_n}}$$ we can simply define $E \subset \prod M_{k_n}$ to be the corresponding extension of $B$ by $\bigoplus M_{k_n}$.  Real rank zero of $E$ follows from the fact that $\prod M_{k_n}$ has real rank zero: every projection in $B$ lifts to a projection in $\prod M_{k_n}$, which necessarily falls in $E$ since $\bigoplus M_{k_n} \triangleleft E$ (cf.\ \cite{brown-pedersen}).

This shows that our initial data exists.

\subsection*{Defining $l(n)$, $r_n$ and $\p_n$}

Now we need some integers, projections and connecting maps.  First let's take natural numbers $m(n)$ which grow so fast that $$\prod_{n\in \mathbb{N}} \frac{m(n) - 1}{m(n)} > 0.$$  Define $l(n) = m(n)k_n$ and identify $M_{l(n)}$ with $M_{m(n)}\otimes M_{k_n}$.  Let $r_n \in M_{m(n)}$ be a projection of rank $m(n) - 1$, but think of it as $r_n \otimes 1_{k_n} \in M_{m(n)}\otimes M_{k_n} = M_{l(n)}$.

At this point, we have taken care of condition (4), so we only have to define the maps $\p_n$ and prove that items (5) - (8) are satisfied.

To ease notation, set $t(n) = l(n)l(n-1)\cdots l(1)$ so that $$M_{t(n)} = M_{l(n)} \otimes \cdots \otimes M_{l(1)}.$$  We define $$\p_n \colon M_{t(n)} \otimes E_n \to M_{l(n+1)} \otimes M_{t(n)}\otimes E_{n+1} = M_{m(n+1)} \otimes M_{k_{n+1}} \otimes M_{t(n)} \otimes E_{n+1}$$ by the formula $$\p_n(T\otimes x) = (r_{n+1}\otimes 1_{k_{n+1}}) \otimes (T \otimes P_{n+1}x) + (r_{n+1}^{\perp}\otimes 1_{k_{n+1}}x) \otimes (T\otimes 1_{E_{n+1}}).$$ Yes, a bit of explanation is in order.

As you probably guessed,  $r_{n+1}^{\perp} = 1_{m(n+1)} - r_{n+1}$.  Recall that $P_{n+1}$ is the central projection in $E_n$ corresponding to the unit of $E_{n+1}$. Hence $x \mapsto P_{n+1}x$ is a well defined unital $*$-homomorphism from $E_n$ to $E_{n+1}$.  The kernel of this morphism is precisely $M_{k_{n+1}} \oplus 0 \oplus 0 \cdots \triangleleft E_n$.  Thus, letting $1_{k_{n+1}}$ denote the central projection in $E_n$ corresponding to its unit, we have a well defined $*$-homomorphism $x \mapsto 1_{k_{n+1}}x$, sending $E_n$ to $M_{k_{n+1}}$.

At this point, condition (5) should be obvious.  The remaining three items require an observation and a rather unpleasant calculation.

The observation is simply that \emph{the map $\p_n \colon M_{t(n)} \otimes E_n \to M_{l(n+1)} \otimes M_{t(n)}\otimes E_{n+1} $ is the canonical inclusion when restricted to $M_{t(n)} \otimes 1_{E_n}$}.  (That is, $\p_n(T\otimes 1) = 1\otimes T \otimes 1$.)  This greatly simplifies the computation of the compositions $\p_{j,n}$.  Indeed, for $T\otimes x \in  M_{t(n)} \otimes E_n$ we have that
\begin{eqnarray*}
\p_{n+2,n}(T\otimes x) &=& (r_{n+2}\otimes 1_{k_{n+2}})\otimes (r_{n+1}\otimes 1_{k_{n+1}}) \otimes (T\otimes P_{n+2}x)\\
& + & (r_{n+2}^{\perp}\otimes 1_{k_{n+2}}x)\otimes (r_{n+1}\otimes 1_{k_{n+1}}) \otimes (T\otimes 1_{E_{n+2}})\\
& + & (1_{l(n+2)}) \otimes   (r_{n+1}^{\perp}\otimes 1_{k_{n+1}}x) \otimes (T\otimes 1_{E_{n+2}}).\\
\end{eqnarray*}

If a TeXnical miracle occurs, the following calculation is without error:
\begin{eqnarray*}
\p_{n+j,n}(T\otimes x) &=& (r_{n+j}\otimes 1_{k_{n+j}}) \otimes \cdots \otimes (r_{n+1}\otimes 1_{k_{n+1}})\otimes (T\otimes P_{n+j} x)\\
&+& (r_{n+j}^{\perp}\otimes 1_{k_{n+j}}x) \otimes (r_{n+j-1}\otimes 1_{k_{n+j-1}})\otimes \cdots \otimes (r_{n+1}\otimes 1_{k_{n+1}})\otimes (T\otimes 1_{E_{n+j}})\\
&+& (1_{l(n+j)})\otimes (r_{n+j-1}^{\perp}\otimes 1_{k_{n+j-1}}x)\otimes \cdots \otimes (r_{n+1}\otimes 1_{k_{n+1}})\otimes (T\otimes 1_{E_{n+j}})\\
&\vdots& \hspace{5cm} \vdots\\
&+& (1_{l(n+j)}) \otimes \cdots \otimes (r_{n+2}^{\perp}\otimes 1_{k_{n+2}}x)\otimes (r_{n+1}\otimes 1_{k_{n+1}})\otimes (T\otimes 1_{E_{n+j}})\\
&+& (1_{l(n+j)}) \otimes \cdots \otimes (1_{l(n+2)}) \otimes (r_{n+1}^{\perp}\otimes 1_{k_{n+1}}x)\otimes (T\otimes 1_{E_{n+j}}).\\
\end{eqnarray*}
If you believe this then conditions (7) and (8) are immediate (thanks to the perpendicular projections appearing in all
the ``scalar" terms).  The last thing to check is (6), but this is also easy because being nonzero means $1_{k_{n+j}}x
\neq 0$ for some $j$.  

This completes the proof of Theorem \ref{thm:main}

\bibliographystyle{amsplain}

\providecommand{\bysame}{\leavevmode\hbox to3em{\hrulefill}\thinspace}

\end{document}